\documentclass[11pt,twoside]{amsart}
\usepackage{amssymb}
\usepackage{mathrsfs}
\usepackage{tikz}
\usepackage{pstricks}

\def\flex#1{\mathrel{\mathop{\kern 0pt\hbox to 
10mm{\rightarrowfill}}\limits_{#1
\rightarrow \infty}}}

\def\Flex#1{\mathrel{\mathop{\kern 0pt\hbox to 
10mm{\rightarrowfill}}\limits_{#1
\rightarrow \infty}}}

\font\bb=msbm10 at 12pt

\def\Z{\hbox{\bb Z}}
\def\Z{\hbox{\bb Z}}

\def\P{\hbox{\bb P}}

\def\T{\hbox{\bb T}}

\def\1{1\hspace{-1.2mm}\mbox{{\normalsize I}}}

\newtheorem{Th}{Theorem}

\newtheorem{Lemme}[Th]{Lemma}

\newtheorem{Prop}[Th]{Proposition}
\newcommand{\R}{\mathbb{R}}

\newcommand{\F}{\mathcal{F}}

\title{ Brownian Motion, Reflection Groups and Tanaka Formula}
\keywords{Reflected Brownian motion; Weyl chamber; reflection groups; longest element; local time; Tanaka formula. \\
{\it AMS classification}: 20F55; 60J65.} 
\begin{document}
\maketitle
\centerline{Nizar Demni\footnote{Universit\'{e} de Rennes I, IRMAR, F-35042 Rennes (nizar.demni@univ-rennes1.fr)} 
and  Dominique L\'{e}pingle\footnote{Universit\'{e} d'Orl\'{e}ans, MAPMO-FDP, F-45067 Orl\'{e}ans (dominique.lepingle@univ-orleans.fr)}} 

\begin{abstract}
In the setting of finite reflection groups, we prove that the projection of a Brownian motion onto a  closed Weyl chamber is another Brownian motion normally reflected on the walls of the chamber. Our proof is probabilistic and the 
decomposition  we obtain may be seen as a 
multidimensional extension of  Tanaka's formula for linear Brownian motion. The paper is closed with a description of the boundary process through the local times of the distances from the initial process to the facets. 
\end{abstract}

\section{Introduction}
The study of stochastic processes in relation to root systems has known a considerable growth during the last decade (see \cite{C2}, \cite{D} and references therein).
In this paper, a reflection group $W$ associated with a root system $R$ is acting on a finite-dimensional Euclidean space $V$, and we project 
a $V$-valued Brownian motion onto the positive closed Weyl chamber $\overline{C}$.
 In fact, each $x \in V$ is conjugated to a unique $y \in \overline{C}$ so that the projection $\pi$ maps vectors in $V$ into their orbits under the action of $W$, namely $\pi:V\rightarrow V/W$. Using the terminology of Dunkl processes, $\pi(X)$ is often called the $W$-invariant or radial  part of the initial Brownian motion $X$ (\cite{C2},\cite{D}), yet should not be confused with the Brownian motion in a Weyl chamber considered in \cite{Bi} which behaves as a Brownian motion conditioned to stay in the open chamber $C$.
The main result of this paper is  the explicit semimartingale  decomposition of the projected process $\pi(X)$. Actually, the latter process is shown to behave as a $V$-valued Brownian motion in the interior of $\overline{C}$ and to reflect normally at the boundary $\partial C$. The proof uses two ingredients of algebraic and probabilistic nature respectively : the expression of $\pi$ as the composition of a finite number of one dimensional projections and Tanaka formula for the absolute value of the linear Brownian motion. Another combination of tools from both reflection groups theory and probability theory leads to a precise description of the boundary process by means of the local time of the distance from $X$ to the facets.

The paper is organized as follows. Section 2 is devoted to basic facts on the theory of reflection groups and to a thorough study of reflections acting on facets of a  root system. An illustration of this study is given in Section 3 where we investigate the evolution of a simple random walk on a triangular lattice in the plane with  the dihedral group $\mathcal{ D}_3$ operating on the lattice. We prove the main result in Section 4. Finally, Section 5 is devoted to the description of the boundary process. 
Note that the problem has previously been tackled in \cite{C1} and taken up again in \cite{D}, however the proof of Theorem 7.1 displayed in \cite{D} p. 216 and quoted from \cite{C1} is not correct. Moreover, our approach is different and more general since we do not assume $\overline{C}$ to be a simplex, an assumption that ceases to hold for some root systems. 

\section{Root systems and reflected facets}
We  give a short introduction to the  theory of reflection groups and refer to \cite{H}  for more details. Let $V$ be an Euclidean space with dimension $N$. For $\alpha \in V$  a unit vector we denote
by $s_{\alpha}$ the orthogonal reflection with respect to the hyperplane $H_{\alpha}:=\{x\in V:\alpha.x=0\}$:
\begin{equation}
s_{\alpha}(x)\,:=\, x\,-\,2\,(\alpha.x)\,\alpha\;.
\end{equation}
A finite subset $R$ of unit vectors in $V$ is called a {\sl reduced root system} if for all $\alpha \in R$
\[
   \begin{array}{l}
   R\cap \R \alpha\,=\,\{\alpha,-\alpha\} \;; \\
   s_{\alpha}(R)\,=\,R \;.
   \end{array}
\]
The group $W\subset O(N)$ which is generated by the reflections $\{s_{\alpha}, \alpha\in R\}$ is called the {\sl reflection group} associated with  $R$. Each hyperplane 
$H_{\beta}:=\{x\in V:\beta.x=0\}$ with $\beta\in V\setminus \cup_{\alpha\in R}H_{\alpha}$ separates the root system $R$ into $R_+$ and $R_-$. Such a set $R_+$ is called a {\sl positive subsystem} and defines the {\sl positive Weyl chamber} $C$ by
\[
 C:= \,\{x\in V:\,\alpha.x>0 \;\; \forall \alpha\in R_+\}\;.
\]
 A subset $S$ of $R_+$ is called {\sl simple} if $S$ is a vector basis for $\textrm{span}(R)$. The elements of $S$ are called {\sl simple roots}. Such a subset exists,  is unique once $R_+$ is fixed and each positive root is a linear combination of simple roots with coefficients being all nonnegative. Moreover (Theorem 1.5 in \cite{H}), $W$ is generated by the reflections $\{s_{\alpha}, \alpha\in S\}$.
The positive Weyl chamber may also be written   
\begin{equation}
  C= \,\{x\in V:\,\alpha.x>0 \;\; \forall \alpha\in S\}\;.
\end{equation}
Its closure is called the {\sl closed Weyl chamber}
\[
  \overline{C}=\,\{x\in V:\,\alpha.x \geq 0 \;\; \forall \alpha\in S\}\;.
\] 

With any partition $R_+=I_+\cup I_0 \cup I_-$ we associate the {\sl facet} 
\[
F=\{x\in V:\, \alpha.x>0 \;\; \forall \alpha \in I_+,\,\beta.x=0 \;\; \forall \beta \in I_0, \,
\gamma.x<0 \;\; \forall \gamma \in I_- \}\;.
\]
We note $(I_+(F), I_0(F),I_-(F))$ the partition corresponding to the facet $F$. We call $\F$ the set of nonempty facets with $Card\{I_0\}=1$. For a facet $F\in {\F}$ with $I_0(F)=\{\beta\}$ , $F\subset H_{\beta}$ and $H_{\beta}$ is called the {\sl support} of $F$. The following result is an easy consequence of the isometry property of $s_{\alpha}$ and of the fact that $s_{\alpha}(R_+\setminus \{\alpha\})= R_+\setminus \{\alpha\}$ (Proposition 1.4 in \cite{H}).
\begin{Lemme} 
For $F\in \F$ and $\alpha\in S$, $s_{\alpha}(F)\in \F$ and
\[
 \begin{array}{lllll}
     I_0(s_{\alpha}(F)) & = &\{\alpha\} &  \quad \mbox{if} & I_0(F)=\{\alpha\} \\
                                  & = & s_{\alpha}(I_0(F)) & \quad \mbox{if} & I_0(F)\neq\{\alpha\} \;.
 \end{array}
\]
For $F\in \F$ and $w\in W$, $w(F)\in \F$ and
\[
  \begin{array}{lllll}
   I_0(w(F)) & = & w(I_0(F))  & \quad \mbox{if} & w(I_0(F))\in R_+ \\
                   & = & -w(I_0(F)) & \quad \mbox{if} & w(I_0(F))\in R_- \;.
  \end{array}
\]
\end{Lemme}

For any $\alpha\in S$ we write $ F_{\alpha}$ for the facet $F$ associated with $I_0=\{\alpha\}, I_-=\emptyset$: 
\[
  F_{\alpha}=\{x\in V:\, \alpha.x=0,\,\beta.x>0 \;\; \forall \beta \in S\setminus \{\alpha\} \}\;.
\]
It is a {\sl face} (\cite{Bo}, p.61) of the chamber $C$ with support the {\sl wall} $H_{\alpha}$. For any simple root $\alpha \in S$ we introduce the mapping $r_{\alpha}$ defined on $V$ by
\[
 \begin{array}{llll}
    r_{\alpha}(x) &  := & x & \mbox{ if }\, \alpha.x\geq 0 \\
                          &  := & s_{\alpha}(x) & \mbox{ if } \, \alpha.x\leq 0 \; .
 \end{array}
\]   
A concise formula is 
\begin{equation}
  r_{\alpha}(x) =\, x+\,2\,(\alpha.x)^- \,\alpha\; . \label{eq:r}
\end{equation}

Let $w\in W$, and let $w=s_{\alpha_l}\cdots s_{\alpha_1}$ be a reduced decomposition of $w$ where  $\alpha_i,i=1,\ldots,l$ are simple roots and $l=l(w)$ is the length of $w$. From Theorem 2.4 in \cite{Bi} and a result of Matsumoto 
(\cite{Bo}, Ch.IV, No.1.5, Prop.5) we know that the operator $r_{\alpha_l}\cdots r_{\alpha_1}$ depends only on $w$ and not on a particular reduced decomposition. 
We denote $r_w$ this operator. Of particular interest is the operator $r_{w_0}$ where $w_0$ denotes the (unique) longest element in $W$: $l(w_0)= \textrm{Card}(R_+)$. The following lemma follows immediately from Corollary 2.9 in \cite{Bi} and   plays a key role in the proof of our main result.

\begin{Lemme}  
If $w_0$ is the longest element in $W$, then $r_{w_0}$ takes values in the closed Weyl chamber $\overline{C}$ and for any $x\in V$, $\pi(x):=r_{w_0}(x)$ is the unique point of the orbit $W.x$ which lies in $\overline{C}$.
\end{Lemme}

The following lemmas shed some light on the action of  $r_{\alpha}$ and $\pi$.

\begin{Lemme}  Let $\alpha \in S$.
\begin{enumerate}
\item
$x\in \cup_{\beta\in R+}H_{\beta} \Longleftrightarrow r_{\alpha}(x) 
\in \cup_{\beta\in R+}H_{\beta}.$
\item
Let  $F\in{\F}$ with support $H_{\beta}$. 
\begin{enumerate}
\item
\begin{itemize}
\item
If $\alpha\in I_+(F)$ , then $r_{\alpha}(F)=F$.
\item
If $\{\alpha\}=I_0(F) $, then $r_{\alpha}(F)=F$ and $\alpha=\beta$.
\item
If $\alpha \in I_-(F)$, then $r_{\alpha}(F)\in {\F}$ and $r_{\alpha}(F)\subset H_{\gamma}$ with $\gamma=s_{\alpha}(\beta)
  \in R_+$.
\end{itemize}
\item
\begin{itemize}
 \item
  If $\alpha\in I_+(F)$ , then $r_{\alpha}^{-1}(F)=F\cup s_{\alpha}(F)$.
\item
If $\{\alpha\}=I_0(F) $, then $r_{\alpha}^{-1}(F)=F$ and $\alpha=\beta$.
\item
If $\alpha \in I_-(F)$, then $r_{\alpha}^{-1}(F)= \emptyset$.
\end{itemize}
\end{enumerate}
\end{enumerate}
 \end{Lemme}
Proof. 
\begin{enumerate}
\item
From the isometry property of $s_{\alpha}\in O(N)$ we derive that for $x\in V$ and $\beta\in R_+$
\[
   \begin{array}{rll}
   \beta . x & = & s_{\alpha}(\beta).s_{\alpha}(x) \\
   s_{\alpha}(\beta).x & = & \beta.s_{\alpha}(x) \;\;\;.
   \end{array}
\]
From Proposition 1.4 in \cite{H} we know that $s_{\alpha}(R_+\setminus \{\alpha\})= R_+\setminus \{\alpha\}$. Both implications
$\Rightarrow$ and $\Leftarrow$ then follow.  
\item
\begin{enumerate}
\item
The proofs of both the first and the second assertions are straightforward. Let $\alpha\in I_-(F)$, $x\in F$ and $\delta \in R_+$. Then,
\[
 \delta.x\,=\,s_{\alpha}(\delta).s_{\alpha}(x)\,=\,s_{\alpha}(\delta).r_{\alpha}(x) \;. 
\]
Therefore,
\[
  \begin{array}{llr}
   \mbox{if} &  \delta=\alpha, & \alpha.r_{\alpha}(x)>0 \\
  \mbox{if} & \delta=\beta, & s_{\alpha}(\beta).r_{\alpha}(x) =0 \\
\mbox{if} & \delta\in I_+(F), &  s_{\alpha}(\delta).r_{\alpha}(x) >0  \\
  \mbox{if} & \delta\in I_-(F)\setminus \{\alpha\}, &  s_{\alpha}(\delta).r_{\alpha}(x) <0  
  \end{array}
\]
and using again Proposition 1.4 in \cite{H} 
we get a new partition $(I_+=s_{\alpha}(I_+(F))\cup \{\alpha\},I_0=\{s_{\alpha}(\beta)\},I_-= s_{\alpha}(I_-(F))\setminus \{\alpha\})$ 
corresponding to
 a new facet $r_{\alpha}(F) \;.$
\item
The  proofs are straightforward or similar to the proofs in (a). $\hfill \square$
\end{enumerate}
\end{enumerate}

\begin{Lemme} 
\begin{enumerate}
\item
 $x\in \cup_{\beta\in R+}H_{\beta} \Longleftrightarrow \pi(x)\in \partial C \;.$
\item
If $F\in \F$, there exists $\alpha\in S$ such that  $\pi(F)= F_{\alpha}$.
\item
 For $\alpha\in S$, $\pi^{-1}(F_{\alpha})=\cup_{F\in{\F}_{\alpha}}F$ where we set 
\[
  {\F}_{\alpha}:=\{F\in {\F}: \pi(F)=F_{\alpha}\} \;.
\] 
\end{enumerate}
\end{Lemme}
Proof. The first assertion is a consequence of (1) of  Lemma 3 and of the property $\pi(V)=\overline{C}$. We also obtain from
(2.a) of Lemma 3 that $\pi(F)\in \F$ 
and it follows there exists $\alpha \in S$ such that $\pi(F)=F_{\alpha}$. For $\alpha\in S$ we deduce from (2.b) of Lemma 3 that 
$\pi^{-1}(F_{\alpha})$ is the union of facets from ${\F}$. Exactly all facets in ${\F}_{\alpha}$ are involved. $\hfill  \square$

\section{Random walk on a triangular lattice} 
Take $V$ to be the Euclidean plane $\R^2$. To each $(i,j)\in \Z^2$ we associate a vertex in the plane with coordinates 
\[
  x_{(i,j)}= i+\frac{1}{2}j \qquad  y_{(i,j)}=\frac{\sqrt{3}}{2}j \;.
\]
The simple random walk $(Z_n)_{n\geq 0}$ on this lattice $\T$  is a Markov chain with uniform transition probability $p((i,j),(k,l))=1/6$ on the six nearest neighbors $(k,l)$ of the vertex $(i,j)$.
We now consider the dihedral group $\mathcal{ D}_3$ consisting of six orthogonal transformations that preserve $\T$. We take $\alpha=(0,1)$ and $\beta=(\sqrt{3}/2,-1/2)$ to be the simple roots, and $\gamma=(\sqrt{3}/2,1/2)$ to be the third positive root. The discrete positive Weyl chamber $C_{dis}$  is given by $\{(i,j)\in \Z^2:i>0,j>0\}$.

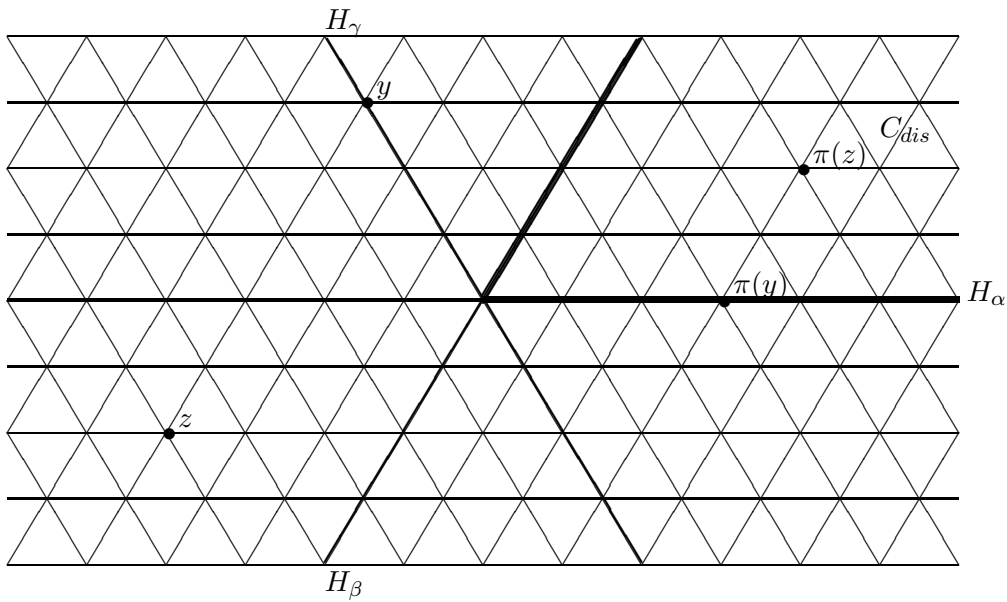
\begin{picture}(400,300)(10,20)
\begin{thinlines}
\multiput(10,90)(30,0){9}{\line(3,5){120}}
\multiput(130,90)(30,0){9}{\line(-3,5){120}}
\multiput(10,90)(0,25){9}{\line(1,0){360}}
\put(280,90){\line(3,5){90}}
\put(310,90){\line(3,5){60}}
\put(340,90){\line(3,5){30}}
\put(370,140){\line(-3,5){90}}
\put(370,190){\line(-3,5){60}}
\put(370,240){\line(-3,5){30}}
\put(100,90){\line(-3,5){90}}
\put(70,90){\line(-3,5){60}}
\put(40,90){\line(-3,5){30}}
\put(10,140){\line(3,5){90}}
\put(10,190){\line(3,5){60}}
\put(10,240){\line(3,5){30}}
\end{thinlines}
\begin{thicklines}
\put(190,190.5){\line(1,0){180}}
\put(190,189.5){\line(1,0){180}}
\put(190.5,191){\line(1,0){179.5}}
\put(190,191){\line(3,5){59}}
\put(189.5,191.5){\line(3,5){58.5}}
\put(190,189){\line(3,5){60}}
\put(10,190.5){\line(1,0){180}}
\put(190,191){\line(-3,5){59}}
\put(130,91){\line(3,5){60}}
\put(250,91){\line(-3,5){60}}
\end{thicklines}
\put(373,190){$H_{\alpha}$}
\put(130,80){$H_{\beta}$}
\put(130,293){$H_{\gamma}$}
\put(340,252){$ C_{dis}$}
\put(71,142){ $z$}
\put(65,137){ $\bullet$}
\put(311,242){ $\pi(z)$}
\put(305,237){ $\bullet$}
\put(146,268){ $y$}
\put(140,262){ $\bullet$}
\put(281,193){ $\pi(y)$}
\put(275,187){ $\bullet$}
\end{picture}

We note that $\pi=r_{\alpha}r_{\beta}r_{\alpha}=r_{\beta}r_{\alpha}r_{\beta}$.
It is easily seen that $(\pi(Z_n))_{n\geq 0}$ is a Markov chain on the set of vertices
$\{(i,j)\in \T:i\geq 0,j\geq 0\}$ with transition probability
\[
 \begin{array}{lll}
  q((i,j),(k,l))= p((i,j),(k,l)) & &\mbox{if} \;  i>0,j>0   \\
  q((i,0),(i,1))=q((i,0),(i-1,1))=2 q((i,0),(i\pm 1,0))   &= \; \frac{1}{3} 
     &   \mbox{if} \; i>0  \\
  q((0,j),(1,j))=q((0,j),(1,j-1))=2 q((0,j),(0,j\pm 1))  &= \;  \frac{1}{3}  
    & \mbox{if}\; j>0  \\
  q((0,0),(0,1))=q((0,0),(1,0))=\frac{1}{2} \;  .
 \end{array}
\]
This is exactly what we should call the simple random walk normally reflected on the walls of the Weyl chamber $C_{dis}$.

The examination of a simple random walk on a square or hexagonal lattice would lead us to the same conclusion.

\section{Brownian motion in a Weyl chamber}
We consider a probability space $(\Omega,\mathscr{A},\P)$ endowed with a right-continuous filtration 
$(\mathscr{A}_t)_{t\geq 0}$ and a $V$-valued Brownian motion $X$ with decomposition
\[
  X_t\,=\, X_0 + B_t \;.
\]
We are interested in the continuous process $\pi(X)$. It can be seen as a $W$-invariant Dunkl process of zero multiplicity function. As such it is a Markov process with values in $\overline{C}$
and semi-group density  (\cite{D}, p.216):
\[
  p_t(x,y)=\frac{1}{c_0t^{N/2}}\exp{\{-\frac{|x|^2+|y|^2}{2t}\}}\sum_{w\in W}\exp{(\frac{1}{t}x.w(y))}
\]
where $c_0$ is a normalizing constant. 
We remark that for any fixed $y\in \overline{C}$ the function $p_t(x,y)$  is a solution to 
the heat equation with Neumann boundary conditions:
\[
  \begin{array}{ll}
   \partial_t p_t= \displaystyle \frac{1}{2} \Delta p_t  &  \mbox{ in } C \\
 \displaystyle \frac{\partial p_t}{\partial n} = 0 & \mbox{ on } \partial C 
  \end{array}
\]
where $\partial p_t/ \partial n$ is any outward normal  derivative of $p_t$. Thus there is no surprise when asserting 
that $\pi(X)$ should be a Brownian motion normally reflected on the boundary of $C$.
 
The aim of this section is to obtain a decomposition of the continuous semimartingale $\pi(X)$.
We first investigate the action of the operator $r_{\alpha}$ on continuous $V$-valued semimartingales. To proceed, recall that if $Z$ is a continuous real semimartingale with local time $L(Z)$ at $0$, then Tanaka formula shows that $Z^-$ is still a semimartingale that decomposes as:
\begin{equation}
Z_t^-=Z_0^--\int_0^t{\bf 1}_{\{Z_s\leq 0\}}dZ_s +\frac{1}{2}L_t(Z)\;  . \label{eq:ta}
\end{equation}

\begin{Lemme} .
Let $Y$ be a continuous $V$-valued semimartingale with Doob decomposition
\[
  Y_t=Y_0+C_t + A_t
\]
where $C$ is a Brownian motion  and A is a continuous process of finite variation, both vanishing at time 0.
Then $Y^{\prime}=r_{\alpha}(Y)$ has a similar decomposition
\[
  Y^{\prime}_t=Y^{\prime}_0+C^{\prime}_t + A^{\prime}_t
\]
where $Y^{\prime}_0=r_{\alpha}(Y_0)$, $C^{\prime}$ is a Brownian motion given by
\[
  C^{\prime}_t=\int_0^tO^{\prime}_s .dC_s
\]
with $O^{\prime}$ a $O(N)$-valued process,
and the finite variation process $A^{\prime}$ is given by
\[
  A^{\prime}_t\,=\,A_t\, -\,2\int_0^t{\bf 1}_{\{\alpha.Y_s\leq 0\}}d(\alpha.A_s)\,\alpha\,+\,L_t(\alpha.Y) \,\alpha \;.
\]
\end{Lemme}
Proof. We combine (\ref{eq:r}) with (\ref{eq:ta}). Then
\[ 
  \begin{array}{ll}
 r_{\alpha}(Y_t) &=\;  Y_t+2 (\alpha.Y_0)^-\,\alpha - 2 \int_0^t{\bf 1}_{\{\alpha .Y_s\leq 0\}} d(\alpha.Y_s) \,\alpha +
                  L_t(\alpha.Y)\,\alpha \\
             & =\;  r_{\alpha}(Y_0) +C_t -2 \int_0^t{\bf 1}_{\{\alpha .Y_s\leq 0\}} d(\alpha.C_s) \,\alpha + A^{\prime}_t 
                          \;.
   \end{array}
\]
The process $O^{\prime}$ defined by
\[
  O^{\prime}_s.u = u-2\,{\bf 1}_{\{\alpha .Y_s\leq 0\}}(\alpha.u)\,\alpha 
\]
 for $u\in V$ satisfies 
\[
 ( O^{\prime}_s.u).(O^{\prime}_s.v)=u.v
\]
for any $u,v\in V$. It follows (\cite{RY}, Ch.IV, Ex.3.22) that
\[
C_t -2 \int_0^t{\bf 1}_{\{\alpha .Y_s\leq 0\}} d(\alpha.C_s) \,\alpha = \int_0^t O^{\prime}_s. dC_s
\]
 defines a $V$-valued Brownian motion. $\hfill  \square$                    

We are now in a position to prove our main result. It states that a
Brownian motion projected onto the closed Weyl chamber is another Brownian motion 
normally reflected on the faces of the chamber.

\begin{Th} .
Let $X_t=X_0+B_t$ be a $V$-valued Brownian motion. Then 
\[
  \pi(X_t)=\pi(X_0)+\int^t_0 O_s.dB_s+\frac{1}{2}\sum_{\alpha\in S}L_t(\alpha.\pi(X))\, \alpha
\]
where $O$ is a $O(N)$-valued process.
\end{Th}  
Proof. \\
 a) Let $\pi=r_{\alpha_l}\cdots r_{\alpha_1}$. We proceed by iteration, setting
$X^0_t:=X_t, B^0_t:=B_t, A_t^0:=0$ and  $X^j_t:=r_{\alpha_j}(X_t^{j-1})$ for $j=1,\ldots,l$. If
we decompose
\[
  X^j_t=X^j_0+B^j_t+A^j_t \;,
\]
then Lemma 5 yields 
\begin{equation}
  \begin{array}{lll}
    X_0^j & = & r_{\alpha_j}(X^{j-1}_0)  \\
   B_t^j & = & \int_0^t O^j_s.dB_s^{j-1} \\
  A_t^j & = & A_t^{j-1} - 2 \int_0^t {\bf 1}_{\{\alpha_j.X^{j-1}_s\leq 0\}} d(\alpha_j.A_s^{j-1}) \, \alpha_j
                  +L_t(\alpha_j.X^{j-1}) \,\alpha_j \;,   
   \end{array}  \label{eq:recu}
\end{equation}
with $O^j$ a $O(N)$-valued process. Setting 
\[
  O := O^l\cdots O^1
\]
we  obtain the Brownian term in the decomposition of $\pi(X)=X^l$. \\
b) From the recurrence formula (\ref{eq:recu})  we deduce that $dA_t^l$
is supported by $\cup_{j=1}^l \{X_t^{j-1}\in H_{\alpha_j}\}$ which is included in
 $ \{X_t\in \cup_{\beta\in R_+} H_{\beta}\}$ by (1) of Lemma 3. 
Moreover, since one-points sets are polar sets for the planar Brownian motion, 
we remark that for any $t>0$, the Brownian motion $X$ does not meet facets with $Card\{I_0\}>1$. Therefore 
\[
  {\bf 1}_{\{X_t\in \cup_{\beta\in R_+} H_{\beta}\}}=\sum_{F\in {\F}} {\bf 1}_{\{X_t\in F\}} \quad\mbox{a.s.}
\] 
and we only
need to study the sequence $(A^j)_{j=1,\ldots,l}$ on the sets $\{X_t\in F\}$ for any $F\in {\F}$. Set $F_0:=F$
and  $F_j:=r_{\alpha_j}(F_{j-1})$ for $j=1,\ldots,l$. Let $H_{\beta_j}$ be the support of the facet $F_j$.
 We will prove by induction that for any $j=0,\ldots,l$
\begin{equation}
{\bf 1}_{\{X_t\in  F\}}dA_t^j = dL^j_t\,\beta_j
\end{equation}
where $L^j$ is an increasing process. This is true for $j=0$ with $L_t^0=0$. Assume this is true for $j-1$ and use (\ref{eq:recu}). Then 
\[
  \begin{array}{l}
  {\bf 1}_{\{X_t\in F\}} dA_t^j \\
  = {\bf 1}_{\{X_t\in F\}} [dL_t^{j-1}\,\beta_{j-1}
    -2\,{\bf 1}_{\{\alpha_j.X^{j-1}_s\leq 0\}}(\alpha_j.\beta_{j-1})\,dL_t^{j-1}\alpha_j +
    dL_t(\alpha_j.X^{j-1})\,\alpha_j] \;.
  \end{array}
\]
We apply Lemma 3.\\
 If $\alpha_j.x>0$ for $x\in F_{j-1}$, then $F_j=F_{j-1}$, $\beta_j=\beta_{j-1}\neq \alpha_j$ and
\[
  {\bf 1 }_{\{X_t\in F\}} dA_t^j = dL_t^{j-1}\,\beta_{j-1}= dL_t^{j-1}\,\beta_j \;.
\]
If $\alpha_j.x=0$ for $x\in F_{j-1}$, then $F_j=F_{j-1}$, $\beta_j=\beta_{j-1}= \alpha_j$ and
\[
  {\bf 1 }_{\{X_t\in F\}} dA_t^j = [-dL^{j-1}_t+ {\bf 1 }_{\{X_t\in F\}} dL_t(\alpha_j.X^{j-1})]\,\beta_j \;.
\]
If $\alpha_j.x<0$ for $x\in F_{j-1}$, then $\alpha_j\neq \beta_{j-1}$, $\beta_j=s_{\alpha_j}(\beta_{j-1})$ and
\[
  {\bf 1 }_{\{X_t\in F\}} dA_t^j =  dL_t^{j-1}\,\beta_j \;.
\]
In the first and third cases  $L^j=L^{j-1}$ is increasing from the induction assumption. 
In the second case  we use Tanaka formula again. As 
\[
\alpha_j.X_t^j
=\alpha_j.r_{\alpha_j}(X^{j-1}_t)\geq 0
\]
we get 
\[
 \begin{array}{lll}
   \alpha_j.X_t^j &= &  (\alpha_j.X_t^j)^+   \\
                  & = & (\alpha_j.X_0^j)^+ +\int_0^t{\bf 1}_{\{\alpha_j.X_s^j>0\}}d( \alpha_j.X_s^j) + \frac{1}{2}\,  L_t(\alpha_j.X^j)
  \end{array}
\]
and therefore 
\[
 \begin{array}{lll}
 dL^j_t & = &{\bf 1 }_{\{X_t\in F\}} d(\alpha_j.A_t^j)\\
                       & = & {\bf 1 }_{\{X_t\in F\}}
                        ({\bf 1}_{\{\alpha_j.X_t^j>0\}}d( \alpha_j.A_t^j) + \frac{1}{2}\, d L_t(\alpha_j.X^j)) \\
                & = &   \frac{1}{2}\,{\bf 1 }_{\{X_t\in F\}}\,  dL_t(\alpha_j.X^j) \;,
  \end{array}
\]
which proves that $L^j$ is an increasing process.
\\
c) We have obtained that  for each $F\in{\F}$ there exists an increasing process $L^l$ such that 
\[
    {\bf 1}_{\{X_t\in F\}} dA_t^l = dL^l_t\,\beta_l \;.
\]
From Lemma 4 we deduce that $F_l=\pi(F)=F_{\alpha}$ for some $\alpha\in S$, which entails $\beta_l=\alpha$.  
 Take $L^{\alpha}$ to be the sum of the $L^l$'s for all $F\in {\F}_{\alpha}$.  Clearly $dL^{\alpha}_t$ is supported by
\[
  \begin{array}{lll}
    \cup_{F\in {\F}_{\alpha}} \{X_t\in F\} & = & \{  X_t\in \cup_{F\in {\F}_{\alpha}}F\} \\
    & = & \{\pi(X_t)\in F_{\alpha}\}
  \end{array}
\]
where the last equality derives from (3) in Lemma 4.
Summarizing we get
\[
  \begin{array}{lll}
    dA_t^l & = & {\bf 1}_{\{X_t\in \cup_{F\in{\F}}\}}dA_t^l \\
    & = & \sum_{\alpha\in S}\sum_{F\in{\F}_{\alpha}}{\bf 1}_{\{X_t\in F\}}dA_t^l \\
   & = & \sum_{\alpha\in S}dL^{\alpha}_t \, \alpha \;. 
  \end{array}
\]
d) It remains to identify the boundary processes $L^{\alpha}$. We use the method in b) again. As 
\[
\alpha.\pi(X_t)\geq 0
\]
we get 
\[
 \begin{array}{lll}
   \alpha.\pi(X_t) &= &  (\alpha.\pi(X_t))^+   \\
                  & = & (\alpha.\pi(X_0))^+ +\int_0^t{\bf 1}_{\{\alpha.\pi(X_s)>0\}}d( \alpha.\pi(X_s)) + \frac{1}{2} L_t(\alpha.\pi(X)) \;.
  \end{array}
\]
On the one hand,
\[
  \begin{array}{lll}
{\bf 1}_{\{\alpha.\pi(X_t)=0\}} d(\alpha.\pi(X_t)) &=& {\bf 1}_{\{\alpha.\pi(X_t)=0\}}(\alpha.(O_s.dB_s)+\sum_{\beta\in S}dL^{\beta}_t\,\alpha.\beta) \\
 & = & dL^{\alpha}_t 
\end{array}
\]
since the set $\{\alpha.\pi(X_t)=0\}$ has zero Lebesgue measure  
and $dL^{\beta}_t$ is supported by $\{\pi(X_t)\in F_{\beta}\}$ for $\beta\neq \alpha$.
On the other hand,
\[
 \begin{array}{lll}
 {\bf 1}_{\{\alpha.\pi(X_t)=0\}} d(\alpha.\pi(X_t))^+ &=& \frac{1}{2}\,dL_t(\alpha.\pi(X))
\end{array}
\]
and we are done. $\hfill  \square $

\section{Complements on the boundary process}
We have already seen in the previous section that on each face $F_{\alpha}$ of the Weyl chamber the boundary process is 
\[
   L_t^{\alpha} = \frac{1}{2}\,dL_t(\alpha.\pi(X))
\]
From  Corollary VI.1.9 in \cite{RY} we know that this is a.s.
\[
\lim_{\varepsilon \downarrow 0}\frac{1}{2\varepsilon} \int_0^t {\bf 1}_{[0,\varepsilon)}(\alpha.\pi(X_s))\,ds
\]
for every $t$. In the sequel, we seek for an expression of $L^{\alpha}$ involving the original process $X$ rather than the reflected one $\pi(X)$. For $x\in V$ and $A\subset V$ we shall use the standard notation 
\[
  d(x,A):=\inf \{|x-y|:\,y\in A\}
\]
and  for $\alpha \in S$ define
\[
  K_{\alpha}:=\pi^{-1}(F_{\alpha})= \cup_{w\in W}w(F_{\alpha}) \;.
\]

\begin{Prop} .
For any $x\in V$ and $\alpha \in S$,
\begin{equation}
\alpha.\pi(x)=d(x,K_{\alpha})\;.   \label{eq:dist}
\end{equation}
\end{Prop}
Proof.
 For any $\varepsilon >0$ there exists $y\in \pi^{-1}(F_{\alpha})$ such that 
\[
  \begin{array}{lll}
  \varepsilon + d(x,K_{\alpha}) & \geq & |x-y| \\
    & \geq & |\pi(x)-\pi(y)| \\
    & \geq & \alpha.(\pi(x)-\pi(y)) \\
    & = & \alpha.\pi(x) \;,
 \end{array}
\]
where we have used the contraction property of $\pi$. In fact  it is easy to check that each $r_{\beta}$ is contracting and so is $\pi$ by iteration, which  proves that $d(x,K_{\alpha})\geq \alpha.\pi(x)$. Conversely, let $x\in V$ and let $w\in W$ be such that 
$w(x)=\pi(x)\in \overline{C}$. Let us consider the facet $w^{-1}(F_{\alpha})$. Using that $\alpha.\beta\leq 0$ for any $\alpha, \beta \in S$ with $\alpha\neq \beta$ (\cite{H}) we  check that
\[
  x-(\alpha .w(x))w^{-1}(\alpha)\in w^{-1}(F_{\alpha}) \subset K_{\alpha}
\]
and therefore
\[
  d(x, K_{\alpha}) \leq \alpha.w(x)=\alpha.\pi(x) \;. 
\]
$\hfill \square$

We obtain a new expression for the boundary process:
\begin{equation}
 L_t^{\alpha} = \lim_{\varepsilon \downarrow 0}\frac{1}{2\varepsilon} \int_0^t {\bf 1}_{\{d(X_s,K_{\alpha})<\varepsilon\}} ds\;.
 \label{eq:bp}
\end{equation}

A little more can be said in a particular case.

\begin{Prop} .
If $\alpha\in S$ is the only simple root in its orbit  $R_{\alpha}:=W.\alpha$, then
\[
  \overline{K_{\alpha}}=\cup_{\gamma\in R_{\alpha}\cap R_+}H_{\gamma} 
\]
and
\[
L_t^{\alpha} = \sum_{\gamma\in R_{\alpha}\cap R_+} L_t(\gamma.X) \;.
\]
\end{Prop}
Proof.
Let $w\in W$ and $F\in {\F}$ satisfy $F=w(F_{\alpha})$. Let $G\in {\F}$ have same support $H_{\gamma}$ ($\gamma\in R_+$) as $F$. There exist $\beta\in S$ and $w^{\prime} \in W$ such that $G=w^{\prime}(F_{\beta})$. Then from Lemma 1
\[
  \gamma = \pm w(\alpha)=\pm w^{\prime}(\beta) \;.
\]
and $\alpha$ and $\beta$ are conjugate. It follows from the hypothesis  that $\beta=\alpha$, and actually the whole hyperplane $H_{\gamma}$ belongs to $\overline{K_{\alpha}}$.
Moreover $\overline{K_{\alpha}}$ is the union of all hyperplanes $H_{\gamma}$ with $\gamma$ a positive root in the orbit of $\alpha$. Now the difference
\[
  \sum_{\gamma\in R_{\alpha}\cap R_+}{\bf 1}_{\{|\gamma .X_s|<\varepsilon\}}
  -{\bf 1}_{\{\min_{\gamma\in R_{\alpha}\cap R_+}|\gamma.X_s|<\varepsilon\}}
\]
is nonnegative and bounded above by a finite sum of terms of the form
 \[
  {\bf 1}_{\{|\gamma_1.X_s|<\varepsilon,|\gamma_2.X_s|<\varepsilon\}} 
\]
where $\gamma_1$ and $\gamma_2$ are independent unit vectors. Let $Y$ be the orthogonal 
projection of $X$ onto the two-dimensional space generated by $\gamma_1$ and $\gamma_2$. Then there exists $c>0$
such that
\[
   {\bf 1}_{\{|\gamma_1.X_s|<\varepsilon,|\gamma_2.X_s|<\varepsilon\}} \leq {\bf 1}_{\{|Y_s|<c\varepsilon\}}
\]
and since the two-dimensional Bessel process $|Y|$ has local time zero at $0$ it follows that a.s.
\[
  \lim_{\varepsilon \downarrow 0}\frac{1}{2\varepsilon} \int_0^t 
  {\bf 1}_{\{|\gamma_1.X_s|<\varepsilon,|\gamma_2.X_s|<\varepsilon\}}ds =0 \;.
\]
Thus
\[
  L_t^{\alpha} =  \lim_{\varepsilon \downarrow 0}\frac{1}{2\varepsilon} \int_0^t 
 \sum_{\gamma\in R_{\alpha}\cap R_+}{\bf 1}_{\{|\gamma .X_s|<\varepsilon\}} ds = 
\sum_{\gamma\in R_{\alpha}\cap R_+} L_t(\gamma.X) \;.      
\]
$\hfill \square$

We now give  two illustrative examples.

 {\bf Dihedral group of order 8}. {\sl  This group consists of the orthogonal transformations which preserve a square centered at the origin in the plane $ V=\R^2$.
There are two simple roots $\alpha= \frac{1}{\sqrt{2}}(e_1-e_2), \beta = e_2$ which lie in two different orbits. The claim of the previous proposition for the simple root $\alpha$ is illustrated by the following picture}:
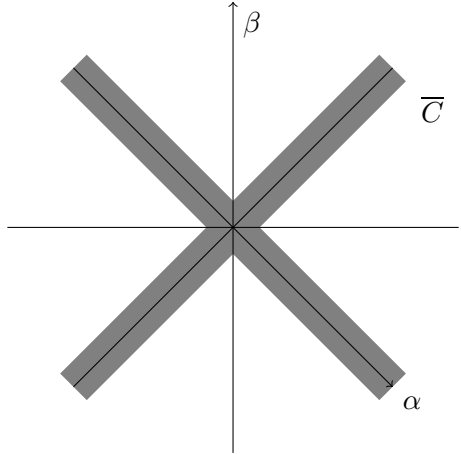
\begin{center}
\begin{tikzpicture}[scale=1]
\draw[line width=5mm, color=gray] (0:0pt) -- (45:3cm) node[at end,below right,color=black]{$\overline{C}$};
\draw[line width=5mm, color=gray] (0:0pt) -- (135:3cm);
\draw[line width=5mm, color=gray] (0:0pt) -- (225:3cm);
\draw[line width=5mm, color=gray] (0:0pt) -- (315:3cm);
\draw[->] (0:0pt) -- (-45:3cm) node[at end,below right]{$\alpha$};
\draw[->] (0:0pt) -- (90:3cm) node[at end,below right]{$\beta$};
\draw[-] (0:0pt) -- (0:3cm);
\draw[-] (0:0pt) -- (135:3cm);
\draw[-] (0:0pt) -- (45:3cm);
\draw[-] (0:0pt) -- (225:3cm);
\draw[-] (0:0pt) -- (-90:3cm);
\draw[-] (0:0pt) -- (180:3cm);
\end{tikzpicture}
\end{center}

{\bf Dihedral group of order 6}. {\sl We come back to the example considered in Section 3 and identify $\R^2$
with the complex plane $\mathbb{C}$. There are six two-dimensional roots represented in complex notations as
\begin{equation}
\{\pm e^{-i\pi/2} e^{im\pi/3}, \, m=1,2,3\}
\end{equation}
and there is only one orbit so that the hypothesis of the proposition breaks down. The reflection group $W = \mathcal{D}_3$ contains three rotations of angles $2m\pi/3, m \in \{1,2,3\}$ and three reflections $\mathbb{C} \ni z \mapsto \overline{z}e^{2im\pi/3}$. When choosing $S = \{\alpha,\beta\} = \{e^{i\pi/2}, e^{-i\pi/6}\}$, $\overline{C}$ is a wedge of angle $\pi/3$ in the positive quadrant of the plane and 
\begin{equation}
 \{ d(x,K_{\alpha})< \varepsilon\}
\end{equation} 
is the hachured part of the following picture:}
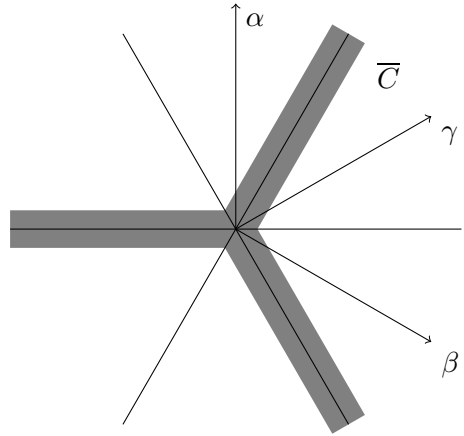
\vspace{0.5cm}
\begin{center}
\begin{tikzpicture}[scale=1]
\draw[line width=5mm, color=gray] (0:0pt) -- (60:3cm) node[at end,below right,color=black]{$\overline{C}$};
\draw[line width=5mm, color=gray] (0:0pt) -- (180:3cm);
\draw[line width=5mm, color=gray] (0:0pt) -- (300:3cm);
\draw[->] (0:0pt) -- (-30:3cm) node[at end,below right]{$\beta$};
\draw[->] (0:0pt) -- (30:3cm) node[at end,below right]{$\gamma$};
\draw[->] (0:0pt) -- (90:3cm) node[at end,below right]{$\alpha$};
\draw[-] (0:0pt) -- (0:3cm);
\draw[-] (0:0pt) -- (120:3cm);
\draw[-] (0:0pt) -- (240:3cm);
\draw[-] (0:0pt) -- (-60:3cm);
\draw[-] (0:0pt) -- (180:3cm);
\draw[-] (0:0pt) -- (60:3cm);
\end{tikzpicture}
\end{center}
It follows that 
\begin{align*}
L_t^{\beta} = \int_0^t {\bf 1}_{\{e^{i\pi/3}. X_s \rangle \geq 0\}}dL_s(\beta.X)&+ \int_0^t {\bf 1}_{\{e^{i\pi}.X_s \geq 0\}} dL_s(\alpha.X)
\\ &+ \int_0^t {\bf 1}_{\{e^{i5\pi/3}. X_s \geq 0\}}dL_s(\gamma.X).
\end{align*}

\vspace{0.7cm}

{\bf Acknowledgements}. We thank J.P. Anker and P. Bougerol  for helpful information  and discussions on reflection groups. We are also grateful to J. Zender for providing us with the Latex code used to draw the last two pictures.

\end{document}